\documentclass[15pt, a4paper]{article}
\usepackage{}
\usepackage{amsfonts}
\usepackage{mathbbold}
\usepackage{bbm}
\usepackage{mathrsfs}
\usepackage{latexsym}
\usepackage{graphicx}
\usepackage{amssymb}

\newtheorem{theorem}{Theorem}[section]
\newtheorem{lemma}[theorem]{Lemma}

\title{{\Large \bf  A note on algebraic connectivity of 2-connected graphs\thanks{Supported by NSFC
(No. 11771376, 11571252, 11601322), Natural Science Foundation of Guangdong Province (2019A1515011031), Foundation of Lingnan Normal
University(ZL1923), ``333" Project of  Jiangsu (2016).}}}
\author{Guanglong Yu$^{a,b}$\thanks{Corresponding authors, E-mail addresses:
yglong01@163.com (G. Yu), rockzhang76@tzc.edu.cn (H. Zhang), wuyarong1@163.com (Y. Wu).}
~ Shuguang Guo$^{b}$ ~ Lin Sun$^{a}$~  Hailiang Zhang$^d$$^{\dag}$ ~ Yarong Wu$^e$$^{\dag}$ ~
\\ ~ \\
{\footnotesize $^a$Department of Mathematics, Lingnan Normal
University,  Zhanjiang, 524048, Guangdong, China}\\ {\footnotesize $^b$Department of Mathematics, Yancheng Teachers
University, Yancheng, 224002, Jiangsu, China}\\
{\footnotesize $^d$Department of Mathematics, Taizhou University, Linhai, Zhejiang, 317000, China}\\
{\footnotesize $^e$ SMU college of art and science, Shanghai maritime
University, Shanghai, 200135, China}}
\date{}

\begin{document}
%\openup 1.0\jot
\maketitle

\begin{abstract}
Among all simple 2-connected graphs, and among all $\theta$-graphs, the graphs with the minimum algebraic connectivity are completely determined, respectively.

\bigskip
\noindent {\bf AMS Classification:} 05C50

\noindent {\bf Keywords:} 2-connected; Algebraic connectivity; $\theta$-graph
\end{abstract}
\baselineskip 18.6pt

\section{Introduction}

\ \ \ \  Algebraic connectivity of a connected graph has been applied extensively to study the real world. It plays an important role for studying the computer science and networks about the robustness, consensus, propagation and so on \cite{AJSU, MWMD, ROFM, SHS, PS}. Consequently, in spectral graph theory, the algebraic connectivity has become an important pole, and a comprehensive researching of the algebraic connectivity of a connected graph has been actively undertaken \cite{MDM, MFL, FKP, GSL, GZY, KRT}.

Given a graph $G$ with vertex set $V(G) = \{v_{0}, v_{1}, \ldots, v_{n-1}\}$ and
edge set $E(G)$, the cardinality $\|V(G)\| =n$ is always called the $order$ of $G$.
Throughout this paper, the graphs considered are connected, undirected and
simple (no loops and no multiple edges).

The $Laplacian$ $matrix$ of graph $G$ is defined to be $L(G) = D(G)-A(G)$, where $A(G)$ is its $adjacency$ $matrix$, $D(G)$ is the $diagonal$ $matrix$ of its degrees.
Because $L(G) = D(G)-A(G)$ is symmetric and positive semidefinite, its
eigenvalues can be assigned as
$\mu_{1}(G) \geq \mu_{2}(G) \geq \cdots \geq \mu_{n-1}(G) \geq \mu_{n}(G) = 0$. The $\mu_{n-1}(G)$ is called $algebraic$ $connectivity$ of $G$.
$\alpha(G)$ is always employed to denote the algebraic connectivity $\mu_{n-1}(G)$. The eigenvectors corresponding to $\alpha(G)$ are always called Fiedler vectors. For a graph $G$, note that all-1-vector $\mathbf{J}=(1$, $1$, $\ldots$, $1)^{T}$ is an eigenvector corresponding to $\mu_{n}(G)$. Hence, it follows that for any $1\leq i\leq n-1$, any eigenvector corresponding to $\mu_{i}(G)$ is orthogonal to $\mathbf{J}$.

It is known from [1] that a graph $G$ is connected if and only if $\alpha(G)>0$.
Moreover, $\alpha(G)$ had been proved to be a good lower bound for the vertex connectivity and edge connectivity of $G$. This
make the algebraic connectivity as a very useful parameter for researching the nice structural properties of a graph or a network.

A lot of nice results about algebraic connectivity of a graph have been shown in \cite{MDM}-\cite{RML} and the references therein. In \cite{GZY}, the authors investigated the algebraic connectivity among all the Hamiltonian graphs of given order and determined the minimum algebraic connectivity. Recently, in \cite{JXLS}, the authors investigated the algebraic connectivity among all the connected graphs with given circumference (the length of the longest cycle) and given order, and determined the minimum algebraic connectivity. With the motivation to investigate the algebraic connectivity about the more general graphs, we consider the 2-connected graphs.

Now, we recall some notions and notations of a graph. In a graph, two paths $P_{1}$ and $P_{2}$ from vertex $u$ to $v$ are called $inner$ $disjoint$ if $V(P_{1})\cap V(P_{2})=\{u, v\}$ (that is, no inner vertex in common). The $local$ $connectivity$ between two distinct vertices $u$ and $v$, denoted by $p(u, v)$, is the maximum number
of pairwise inner disjoint paths from $u$ to $v$. A nontrivial graph $G$ is 2-connected if $p(v_{1}, v_{2})\geq 2$ for its any two distinct vertices
$v_{1}$ and $v_{2}$. A $\theta$-graph is a 2-connected graph which consists of three pairwise inner disjoint paths with common initial and terminal vertices.

\setlength{\unitlength}{0.6pt}
\begin{center}
\begin{picture}(368,444)
\put(35,390){\circle*{4}}
\put(95,358){\circle*{4}}
\qbezier(35,390)(65,374)(95,358)
\put(95,423){\circle*{4}}
\qbezier(35,390)(65,407)(95,423)
\put(152,423){\circle*{4}}
\qbezier(95,423)(123,423)(152,423)
\put(152,358){\circle*{4}}
\qbezier(95,358)(123,358)(152,358)
\qbezier(152,423)(152,391)(152,358)
\qbezier(95,423)(95,391)(95,358)
\put(190,419){\circle*{4}}
\put(206,419){\circle*{4}}
\put(221,419){\circle*{4}}
\put(191,362){\circle*{4}}
\put(207,362){\circle*{4}}
\put(222,362){\circle*{4}}
\put(254,423){\circle*{4}}
\put(254,361){\circle*{4}}
\qbezier(254,423)(254,392)(254,361)
\put(306,423){\circle*{4}}
\qbezier(254,423)(280,423)(306,423)
\put(306,361){\circle*{4}}
\qbezier(254,361)(280,361)(306,361)
\qbezier(306,423)(306,392)(306,361)
\put(15,388){$v_{0}$}
\put(90,344){$v_{1}$}
\put(148,344){$v_{2}$}
\put(296,345){$v_{\frac{n-1}{2}}$}
\put(298,435){$v_{\frac{n+1}{2}}$}
\put(145,430){$v_{n-2}$}
\put(82,430){$v_{n-1}$}
\put(173,314){$H_{1}$}
\put(23,242){\circle*{4}}
\put(79,277){\circle*{4}}
\qbezier(23,242)(51,260)(79,277)
\put(79,209){\circle*{4}}
\qbezier(23,242)(51,226)(79,209)
\qbezier(79,277)(79,243)(79,209)
\put(135,277){\circle*{4}}
\qbezier(79,277)(107,277)(135,277)
\put(135,209){\circle*{4}}
\qbezier(79,209)(107,209)(135,209)
\qbezier(135,277)(135,243)(135,209)
\put(171,272){\circle*{4}}
\put(187,272){\circle*{4}}
\put(202,272){\circle*{4}}
\put(172,215){\circle*{4}}
\put(188,215){\circle*{4}}
\put(203,215){\circle*{4}}
\put(238,276){\circle*{4}}
\put(238,211){\circle*{4}}
\qbezier(238,276)(238,244)(238,211)
\put(290,276){\circle*{4}}
\qbezier(238,276)(264,276)(290,276)
\put(290,211){\circle*{4}}
\qbezier(238,211)(264,211)(290,211)
\qbezier(290,276)(290,244)(290,211)
\put(347,243){\circle*{4}}
\qbezier(290,276)(318,260)(347,243)
\qbezier(290,211)(318,227)(347,243)
\put(5,241){$v_{0}$}
\put(73,196){$v_{1}$}
\put(129,196){$v_{2}$}
\put(350,241){$v_{\frac{n}{2}}$}
\put(277,288){$v_{\frac{n+2}{2}}$}
\put(280,199){$v_{\frac{n-2}{2}}$}
\put(68,286){$v_{n-1}$}
\put(174,171){$H_{2}$}
\put(21,96){\circle*{4}}
\put(75,132){\circle*{4}}
\qbezier(21,96)(48,114)(75,132)
\put(66,62){\circle*{4}}
\qbezier(21,96)(43,79)(66,62)
\put(136,132){\circle*{4}}
\qbezier(75,132)(105,132)(136,132)
\put(135,62){\circle*{4}}
\qbezier(66,62)(100,62)(135,62)
\put(200,132){\circle*{4}}
\qbezier(136,132)(168,132)(200,132)
\qbezier(136,132)(101,97)(66,62)
\qbezier(200,132)(167,97)(135,62)
\put(232,126){\circle*{4}}
\put(248,126){\circle*{4}}
\put(263,126){\circle*{4}}
\put(184,69){\circle*{4}}
\put(200,69){\circle*{4}}
\put(215,69){\circle*{4}}
\put(301,130){\circle*{4}}
\put(257,63){\circle*{4}}
\qbezier(301,130)(279,97)(257,63)
\put(358,97){\circle*{4}}
\qbezier(301,130)(329,114)(358,97)
\put(315,63){\circle*{4}}
\qbezier(257,63)(286,63)(315,63)
\qbezier(315,63)(336,80)(358,97)
\put(3,95){$v_{0}$}
\put(58,50){$v_{1}$}
\put(129,50){$v_{2}$}
\put(362,96){$v_{\frac{n}{2}}$}
\put(64,139){$v_{n-1}$}
\put(172,22){$H_{3}$}
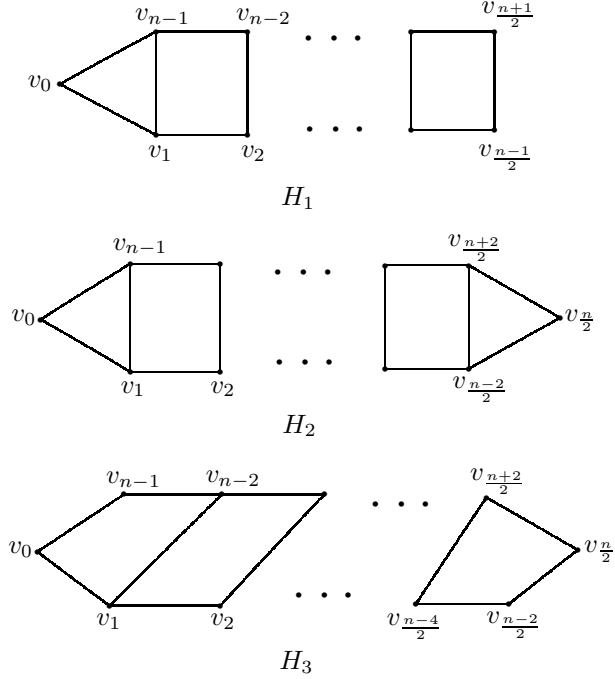
\put(129,-9){Fig 1.1. $H_{1}$, $H_{2}$, $H_{3}$}
\put(239,50){$v_{\frac{n-4}{2}}$}
\put(301,51){$v_{\frac{n-2}{2}}$}
\put(291,141){$v_{\frac{n+2}{2}}$}
\put(126,140){$v_{n-2}$}
\end{picture}
\end{center}

Given a graph $G$, let $G +uv$ denote the graph
obtained from $G$ by adding a new edge $uv\notin E(G)$ between two nonadjacent vertices $u, v$ in $G$; let $G-uv$ denote the graph
obtained from $G$ by deleting an edge $uv\in E(G)$; for a graph $K$ with $\|E(K)\|\geq 2$, $V(K)\subseteq  V(G)$ and $E(K)\nsubseteq  E(G)$, let $G +K$ denote the graph
obtained from $G$ with new edge set $E(G+K)=E(G)\cup E(K)$, where $E(K)\cap  E(G)\neq \emptyset$ possibly. Denote by $C_{n} = v_{0}v_{ 1} v_{ 2}\cdots v_{n-1}v_{0}$ a cycle with $n$ vertices. Let $H_{1} (i_{ 1}, \ldots, i_{ k})= C_{ n} +v_{ i_{ 1}} v_{ n-i_{1}} +\cdots+v_{ i_{ k}} v_{ n-i_{ k}} $, where $n \geq 5$ is odd, $1 \leq k \leq \frac{n -3}{2}$
and $1 \leq i_{ 1} \leq \cdots \leq i_{ k} \leq \frac{n -3}{2}$; $H_{1}=H_{1} (1, 2, \ldots, \frac{n -3}{2})$ (see Fig. 1.1). Let $H_{2} (i_{ 1}, \ldots, i_{ k})= C_{ n} +v_{ i_{ 1}} v_{ n-i_{1}} +\cdots+v_{ i_{ k}} v_{ n-i_{ k}} $, where $n \geq 4$ is even, $1 \leq k \leq \frac{n-2}{2}$
and $1 \leq i_{ 1} \leq \cdots \leq i_{ k} \leq \frac{n-2}{2}$; $H_{2}=H_{2} (1, 2, \ldots, \frac{n-2}{2})$ (see Fig. 1.1). Let $H_{3} (i_{ 1}, \ldots, i_{ k})= C_{ n} +v_{ i_{ 1}} v_{ n-i_{1}-1} +\cdots+v_{ i_{ k}} v_{ n-i_{ k}-1} $, where $n \geq 6$ is even, $1 \leq k \leq \frac{n -4}{2}$
and $1 \leq i_{ 1} \leq \cdots \leq i_{ k} \leq \frac{n -4}{2}$; $H_{3}=H_{3} (1, 2, \ldots, \frac{n -4}{2})$ (see Fig. 1.1).

In this paper, for determining the minimum algebraic connectivity among all the simple 2-connected graphs and among all $\theta$-graphs, we get the following results:

\begin{theorem}\label{le03.02} %------
Let $G$ be a 2-connected graph of order $n$. Then
$\alpha(G)\geq \alpha(C_{n})$. Moreover,

(i) if the order $n$ is odd and the equality holds, then $G \cong C_{ n}$ or $G \cong H_{1} (i_{ 1}, \ldots, i_{ k})$ for some $i_{ 1}$, $\ldots$, $i_{ k}$;

(ii) if the order $n$ is even and the equality holds, then $G \cong C_{ n}$, $G \cong H_{2} (i_{ 1}, \ldots, i_{ k})$ for some $i_{ 1}$, $\ldots$, $i_{ k}$, or $G \cong H_{3} (j_{ 1}, \ldots, j_{s})$ for some $j_{ 1}$, $\ldots$, $j_{ s}$.
\end{theorem}

\begin{theorem}\label{le03.03} %------
Let $G$ be a $\theta$-graph of order $n\geq 4$. Then $\alpha(G)\geq\alpha(C_{n})$. Moreover,

(i) if the order $n$ is odd and the equality holds, then $G \cong H_{1} (i)$ for some $1\leq i\leq (n -3)/2$;

(ii) if the order $n$ is even and the equality holds, then $G \cong H_{2} (i)$ for some $1\leq i\leq \frac{n-2}{2}$, or $G \cong H_{3} (i)$ for some $1\leq i\leq \frac{n-4}{2}$.
\end{theorem}

\section{Proof of main results}

\begin{lemma}{\bf \cite{FRGR}}\label{le02.01} %------
Let $G$ be a connected graph with $n$ vertices $v_{0}$, $v_{1}$, $\ldots$, $v_{n-1}$. Then
$$\displaystyle \alpha(G)=\min_{(X^{T}, \mathbf{J})=0}\frac{\displaystyle\sum_{v_{i}v_{k}\in E(G)}(x(v_{i})-x(v_{k}))^{2}}{\displaystyle\sum_{v_{j}\in V(G)}x^{2}(v_{j})},$$ where $X=(x(v_{0})$, $x(v_{1})$, $\ldots$, $x(v_{n-1}))^{T}$ is a vector in which $x(v_{i})$ corresponds to vertex $v_{i}$ for $0\leq i\leq n-1$, and where $(X^{T}, \mathbf{J})$ is the usual inner product.
\end{lemma}

\begin{lemma}{\bf \cite{GZY}}\label{le02,04} %------
Let $G$ be a Hamiltonian graph with $n$ vertices.
Then
$\alpha(G)\geq \alpha(C_{n})$. Moreover,

(i) if the order $n$ is odd and the equality holds, then $G \cong C_{ n}$ or $G \cong H_{1} (i_{ 1}, \ldots, i_{ k})$;

(ii) if the order $n$ is even and the equality holds, then $G \cong C_{ n}$, $G \cong H_{2} (i_{ 1}, \ldots, i_{ k})$ for some $i_{ 1}$, $\ldots$, $i_{ k}$, or $G \cong H_{3} (j_{ 1}, \ldots, j_{s})$ for some $j_{ 1}$, $\ldots$, $j_{ s}$.

\end{lemma}

{\bf Proof of Theorem \ref{le03.02}}
Let $V(G) = \{v_{0}, v_{1}, \ldots, v_{n-1}\}$, $X=(x(v_{0})$, $x(v_{1})$, $\ldots$, $x(v_{n-1}))^{T}$ be
a Fiedler vector of $G$ where $x(v_{i})$ corresponds to vertex $v_{i}$ for $0\leq i\leq n-1$, $x(v_{\tau})=\max\{x(v_{i})\mid v_{i}\in V(G)\}$ and $x(v_{0})=\min\{x(v_{i})\mid v_{i}\in V(G)\}$. It is known that $x(v_{0})<x(v_{\tau})$ because $X$ is orthogonal to $\mathbf{J}$. Note that $G$ is 2-connected. we get that there are two inner disjoint paths $P_{1}$ and $P_{2}$ from $v_{0}$ to $v_{\tau}$. Then $C=P_{1}\cup P_{2}$ is a cycle. If  $\|V(C)\|= n$, then $G$ is a Hamiltonian graph, and then the results of Theorem \ref{le03.02} follows from Lemma \ref{le02,04} directly. Next, we suppose $\|V(C)\|< n$.

Suppose that $P_{1}=v_{0}v_{i_{1}}v_{i_{2}}\cdots v_{i_{s}}v_{\tau}$. Here, for convenience, we let $i_{0}=0$, $i_{s+1}=\tau$. Note that $x(v_{0})<x(v_{\tau})$. Then there is $0\leq w\leq s$ that $x(v_{i_{w}})\neq x(v_{i_{w+1}})$. For this $w$, if $x(v_{i_{w}})< x(v_{i_{w+1}})$ and $v_{a_{1}}$, $v_{a_{2}}$, $\ldots$, $v_{a_{g}}$ are all the vertices in $V(G)\setminus V(C)$ that $x(v_{i_{w}})< x(v_{a_{1}})\leq x(v_{a_{2}})\leq\cdots\leq x(v_{a_{g}})\leq x(v_{i_{w+1}})$, then we let path $P_{w, w+1}=v_{i_{w}}v_{a_{1}} v_{a_{2}}\cdots v_{a_{g}}v_{i_{w+1}}$; if $x(v_{i_{w}})> x(v_{i_{w+1}})$ and $v_{a_{1}}$, $v_{a_{2}}$, $\ldots$, $v_{a_{g}}$ are all the vertices in $V(G)\setminus V(C)$ that $x(v_{i_{w+1}})< x(v_{a_{1}})\leq x(v_{a_{2}})\leq\cdots\leq x(v_{a_{g}})\leq x(v_{i_{w}})$, then we let $P_{w, w+1}=v_{i_{w+1}}v_{a_{1}} v_{a_{2}}\cdots v_{a_{g}}v_{i_{w}}$. Here, $g=0$ means that for any vertex $v\in (V(G)\setminus V(C))$, neither $x(v_{i_{w}})< x(v)\leq x(v_{i_{w+1}})$ holds, nor $x(v_{i_{w+1}})< x(v)\leq x(v_{i_{w}})$ holds. Therefore, if $g=0$, then $\|V(P_{w, w+1})\|=2$, $P_{w, w+1}=v_{i_{w}}v_{i_{w+1}}$ or $P_{w, w+1}=v_{i_{w+1}}v_{i_{w}}$. Suppose for $1\leq \epsilon\leq d$, $P_{w_{\epsilon},w_{\epsilon}+1}$ are the paths that $\|V(P_{w_{\epsilon},w_{\epsilon}+1})\|> 2$, and $V(G)\setminus V(C)=\cup^{d}_{\epsilon=1}(V(P_{w_{\epsilon},w_{\epsilon}+1})\setminus \{v_{i_{w_{\epsilon}}}, v_{i_{w_{\epsilon}+1}}\})$, where $0\leq w_{1}\leq w_{2}\leq \cdots \leq w_{d}\leq s$, $(V(P_{w_{j},w_{j}+1})\setminus \{v_{i_{w_{j}}}, v_{i_{w_{j}+1}}\})\cap (V(P_{w_{k},w_{k}+1})\setminus \{v_{i_{w_{k}}}, v_{i_{w_{k}+1}}\})=\emptyset$ for $1\leq j\neq k\leq d$.

Let $G^{'}=G-\sum^{d}_{\epsilon=1} v_{i_{w_{\epsilon}}}v_{i_{w_{\epsilon}+1}}+\sum^{d}_{\epsilon=1} P_{w_{\epsilon},w_{\epsilon}+1}$. Note that for any three real numbers $h\leq l\leq q$, we have $(q-h)^{2}\geq (q-l)^{2}+(l-h)^{2}$. And from this, for $h\leq l_{1}\leq l_{2}\leq \cdots \leq l_{z}\leq q$, we have $(q-h)^{2}\geq (q-l_{z})^{2}+(l_{1}-h)^{2}+\sum^{z-1}_{j=1}(l_{j+1}-l_{j})^{2}$ if $z\geq 1$. Thus for $1\leq \epsilon\leq d$, if $x(v_{i_{w_{\epsilon}}})< x(v_{i_{w_{\epsilon}+1}})$, then for $P_{w_{\epsilon},w_{\epsilon}+1}=v_{i_{w_{\epsilon}}}v_{\epsilon_{1}}v_{\epsilon_{2}}\cdots v_{\epsilon_{g(\epsilon)}}v_{i_{w_{\epsilon}+1}}$ (where $g(\epsilon)\geq 1$), it follows that $$(x(v_{\epsilon_{1}})-x(v_{i_{w_{\epsilon}}}))^{2}+\sum^{g(\epsilon)-1}_{j=1}(x(v_{\epsilon_{j+1}})-x(v_{\epsilon_{j}}))^{2}+(x(v_{i_{w_{\epsilon}+1}})-x(v_{\epsilon_{g(\epsilon)}}))^{2}\leq (x(v_{i_{w_{\epsilon}+1}})-x(v_{i_{w_{\epsilon}}}))^{2};$$ if $x(v_{i_{w_{\epsilon}}})> x(v_{i_{w_{\epsilon}+1}})$, then for $P_{w_{\epsilon},w_{\epsilon}+1}=v_{i_{w_{\epsilon}+1}}v_{\epsilon_{1}}v_{\epsilon_{2}}\cdots v_{\epsilon_{g(\epsilon)}}v_{i_{w_{\epsilon}}}$, it follows that $$(x(v_{\epsilon_{1}})-x(v_{i_{w_{\epsilon}+1}}))^{2}+\sum^{g(\epsilon)-1}_{j=1}(x(v_{\epsilon_{j+1}})-x(v_{\epsilon_{j}}))^{2}+(x(v_{i_{w_{\epsilon}}})-x(v_{\epsilon_{g(\epsilon)}}))^{2}\leq (x(v_{i_{w}})-x(v_{i_{w+1}}))^{2}.$$ Consequently, it follows that
$X^{T}L(G^{'})X\leq X^{T}L(C)X\leq X^{T}L(G)X$. Then by Lemma \ref{le02.01}, we get that $\alpha(G^{'})\leq \frac{X^{T}L(G^{'})X}{X^{T}X}\leq \frac{X^{T}L(G)X}{X^{T}X}=\alpha(G)$.

Suppose $\alpha(G^{'})= \alpha(G)$. Then $\alpha(G^{'})= \frac{X^{T}L(G^{'})X}{X^{T}X}= \frac{X^{T}L(G)X}{X^{T}X}=\alpha(G)$. Thus $X$ is also a Fiedler vector of $G^{'}$. Without loss of generality, suppose that $P_{w_{d}, w_{d}+1}=v_{i_{w_{d}}}v_{d_{1}}v_{d_{2}}\cdots v_{d_{g(d)}}v_{i_{w_{d}+1}}$ where $g(d)\geq 1$, $x(v_{i_{w_{d}+1}})>x(v_{i_{w_{d}}})$. Note that $$\alpha(G^{'})x(v_{i_{w_{d}+1}})=(L(G^{'})X)_{v_{i_{w_{d}+1}}}=(L(G)X)_{v_{i_{w_{d}+1}}}-x(v_{i_{w_{d}+1}})+x(v_{i_{w_{d}}})$$
$$\hspace{5.3cm}=\alpha(G)x(v_{i_{w_{d}+1}})-x(v_{i_{w_{d}+1}})+x(v_{i_{w_{d}}}),$$ if $v_{d_{g(d)}}$ is adjacent to $v_{i_{w_{d}+1}}$ in $G$ (where $(L(G^{'})X)_{v_{i_{w_{d}+1}}}$ is the element corresponds to vertex $v_{i_{w_{d}+1}}$ in $L(G^{'})X$); $$\alpha(G^{'})x(v_{i_{w_{d}+1}})=(L(G^{'})X)_{v_{i_{w_{d}+1}}}=(L(G)X)_{v_{i_{w_{d}+1}}}-x(v_{d_{g(d)}})+x(v_{i_{w_{d}}})$$
$$\hspace{5.3cm}=\alpha(G)x(v_{i_{w_{d}+1}})-x(v_{d_{g(d)}})+x(v_{i_{w_{d}}})$$ if $v_{d_{g(d)}}$ is not adjacent to $v_{i_{w_{d}+1}}$ in $G$. Note that $x(v_{i_{w_{d}+1}})>x(v_{i_{w_{d}}})$ and $x(v_{d_{g(d)}})>x(v_{i_{w_{d}}})$. Thus it follows that $\alpha(G^{'})\neq \alpha(G)$, which contradicts the supposition that $\alpha(G^{'})= \alpha(G)$. Then we get that $\alpha(G^{'})< \alpha(G)$.

From the above discussion, combining Lemma \ref{le02,04}, for the case that $\|V(C)\|< n$, we get that the results in Theorem \ref{le03.02} hold.
The proof is completed. \ \ \ \ \ $\Box$

Theorem \ref{le03.03} follows from Theorem \ref{le03.02} as a corollary.

\small {

}

\end{document}